\documentclass[12pt]{article}
\newtheorem{theorem}{Theorem}[section]
\newtheorem{lemma}[theorem]{Lemma}

\renewcommand{\sp}{\vspace{1ex}}

\newcommand{\pl}{\partial}

\newcommand{\ov}{\over}
\newcommand{\ra}{\rightarrow}

\newcommand{\proof}{{\em Proof. }}
\newcommand{\be}{\begin{equation}}
\newcommand{\ee}{\end{equation}}
\newcommand{\bq}{\begin{eqnarray}}
\newcommand{\eq}{\end{eqnarray}}

\newcommand{\ba}{\begin{array}}
\newcommand{\ea}{\end{array}}

\renewcommand{\Re}{\mbox{\rm Re\,}}

\newcommand{\inv}{^{-1}}
\newcommand{\iy}{\infty}
\newcommand{\al}{\alpha}

\newcommand{\tr}{\mbox{\rm tr\,}}

\newcommand{\nc}{\newcommand}   \nc{\cd}{\cdots}
\renewcommand{\sp}{\vskip1ex} \nc{\noi}{\noindent}
\nc{\onetwo}[2]{\left(\begin{array}{cc}#1&#2\end{array}\right)}
\nc{\twotwo}[4]{\left(\begin{array}{cc}#1&#2\\&\\#3&#4\end{array}\right)}
\nc{\twoone}[2]{\left(\begin{array}{c}#1\\\\#2\end{array}\right)}
\renewcommand{\l}{L_n^{\nu}}
\nc{\J}{J_{\nu}}
\nc{\Jm}{J_{-\nu}} \nc{\st}[1]{2\sqrt{ N#1}} \nc{\stp}[1]{2\sqrt{ N'#1}}

\begin{document}

\title{\bf Determinants of Hankel Matrices}
\author{Estelle L. Basor\thanks{Supported in part by NSF Grants
DMS-9623278 and DMS-9970879.}\\
               \emph {Department of Mathematics}\\
               \emph {California Polytechnic State University}\\
               \emph {San Luis Obispo, CA 93407, USA}
        \and
     Yang Chen\\
               \emph {Department of Mathematics}\\
                   \emph {Imperial college}\\
                   \emph {180 Queen's Gate, London, SW7 2BZ, UK}
        \and
     Harold Widom\thanks{Supported in part by NSF Grant DMS-9732687}\\
               \emph {Department of Mathematics}\\
                   \emph {University of California}\\
                   \emph {Santa Cruz, CA 96054, USA}}
\date{June 8, 2000}

\maketitle

\begin{abstract}
The purpose of this paper is to compute asymptotically Hankel
determinants
for weights that are supported in a semi-infinite interval.
The main idea is to reduce the problem to determinants of other
operators whose determinant asymptotics are well known.
\end{abstract}

%%%%%%%%%%%%%%%%%%%%%%%%%%%%%%%%%%%%%%%%%%%%%%%%%%%%%%%%%%%%%%%%%%%

\section{Introduction}

The main purpose of this paper is to compute asymptotically the
determinants of the finite matrices $H_{n}(u)$ defined by
\[\det(a_{i+j})_{i,j=0}^{n-1}\] where
\[a_{i+j}=\int_{0}^{\iy}x^{i+j}u(x)\,\,dx\]
with suitable conditions on the weight function $u(x).$
These determinant entries depend only the sum $i+j$ and are hence
classified
as Hankel determinants.

Hankel determinants such as these were considered by Szeg\"{o} in
\cite{Sz2} and also by Hirshmann in \cite{H}, but in both cases for
finite intervals.
These determinants are important in random matrix theory and its
applications.

Our main result is as follows. Suppose we replace $u(x)$ by a function
given in the form $w(x)U(x)$ where $w(x)$ is the weight $e^{-x}x^{\nu}$
with $\nu\ge-1/2$.
Then if $U$ is nowhere zero and $U-1$ is a Schwartz function (a
condition which
can be considerably relaxed) the determinants are given
asymptotically as $n\ra \iy$ by
\bq
\det (H_{n}(u)) =\exp\{c_{1}n^{2}\log n +c_{2}n^{2} +c_{3}n\log n
+c_{4} n +c_{5}n^{1/2} +c_{6}\log n + c_{7}+o(1)\} \label{asyformula}
\eq
where
$$c_{1}= 1,\,\,\,\,\,
 c_{2}=-3/2, \,\,\,\,\,
 c_{3}=\nu,$$
$$c_{4}=-\nu +\log 2\pi,\,\,\,\, c_{5}=  \frac{2}{\pi}
\int_{0}^{\iy}\log(U(x^{2})) dx$$
$$c_{6}=\nu^{2}/2 -1/6, \,\,\,\,\,
 c_{7}=4/3\log G(1/2) +(1/3 +\nu/2)\log \pi+ (\nu/2 -1/18)\log 2$$
 $$ -\log G(1+\nu)
 -\nu/2\log U(0) +
\frac{1}{2\pi^{2}}\int_{0}^{\iy}xS(x)^{2}dx,$$
$G$ is the Barnes G-function, and $S(x) =\int_{0}^{\iy}\cos(xy) \log
U(y^{2})dy.$

The idea behind the proof is to replace the matrix $H_{n}(u)$ with one
whose $i,j^{th}$ entry is given by
\[\int_{0}^{\iy}P_{i}(x)P_{j}(x)u(x)\,\,dx\] where the $P_{i}$'s are
orthogonal (Laguerre) polynomials with respect to the weight
$e^{-x}x^{\nu}.$

These new
determinants can then be evaluated using the ideas of the ``linear
statistic'' method in random matrix theory. More precisely, the above
determinant can be replaced by
\[\det (I+K_{n})\] where $K_{n}$ is an integral operator whose
kernel is given by
\[\sum_{i=0}^{n-1}L_{i}^{\nu}(x)L_{i}^{\nu}(y)(U-1)(y)\] with
$L_{i}^{\nu}$
defined as the $i^{th}$ Laguerre function of order $\nu.$ The main
computation in the paper is to approximate the above kernel
with a different kernel which involves Bessel functions. This new
kernel was fortunately already considered in \cite{B}. There,
asymptotics
for certain integral operators were computed and these results
are then applied to give the result of formula (\ref{asyformula}).
The  kernels considered in \cite{B} arises in random matrix theory in
the
``hard-edge'' scaling for ensembles of positive Hermitian matrices.
Details about the random matrix connections can be found in \cite{B}.

The formula in (1) was stated earlier in \cite{BCW} where a heuristic
argument using the Coulomb gas approach was used to derive the same
formula. The Coulomb gas approach was also used in \cite{BCW} to extend
to the case
where the interval of integration is the entire real line

\section{Preliminaries}
We first show how to replace the powers in the Hankel determinants
with the orthogonal polynomials.
Let \[{\cal H}_{n}(u) =
\left(\int_{0}^{\iy}P_{i}(x)P_{j}(x)u(x)\,\,dx\right)\biggl|_{i,j=0}^{n-1},\]
and write \[P_{i}(x) = \sum_{k \leq i}a_{ik}x^{k}.\]
 \begin{lemma} Let $P_{i}, {\cal H}_{n}(u), H_{n}(u)$ be defined as
above
 and let us assume that all moment integrals exist.
 Then \[\det({\cal H}_{n}(u)) = A_{n}\det(H_{n}(u)),\] where \[A_{n} =
 (\prod_{i=0}^{n-1}a_{ii})^{2}.\]
 \end{lemma}
 \proof
 We have \[\det({\cal H}_{n}(u) )= \det (\int_{0}^{\iy}(\sum_{m \leq
j}a_{jm}x^{m})
 (\sum_{l \leq k }a_{kl}x^{l}) u(x) \,\,dx).\] The $j,k$ entry of the
 matrix is given by
 \[\sum_{m\leq j}\sum_{l \leq
k}a_{jm}a_{kl}\int_{0}^{\iy}x^{m+l}u(x)\,\,dx\]
 \[ =\sum_{m\leq j}\sum_{l \leq k}a_{jm}\int_{0}^{\iy}x^{m+l}u(x)\,\,dx
\,\,a_{kl}.
 \]
Form this it follows that ${\cal H}_{n}(u)$ is a product of three
matrices $T H_{n}(u) T^{t}.$ The matrix $T$ is lower triangular with
entries $a_{jm}, m \leq j.$ It is easy to see from this that the lemma
follows.

The next step is to evaluate the term
\[A_{n} =
 (\prod_{i=0}^{n-1}a_{ii})^{2},\]
 which is of course straight forward. We normalize the
 polynomials so that they are orthonormal. Then it is well known
 \cite{Sz} that \[a_{ii}^{2} =
 \frac{1}{\Gamma(1+i+\nu)\Gamma(1+i)}.\]
 The following
 lemma computes the product asymptotically using the Barnes G-function.
 This function is defined by \cite{Bar,WW}
\bq
 G(1+z) &=&(2\pi)^{z/2}e^{-(z+1)z/2-\gamma_{E}z^{2}/2}
 \prod_{k=1}^\iy\left((1+z/k)^{k}e^{-z+z^{2}/2k}\right)
\eq
with $\gamma_{E}$ being Euler's constant.
\begin{lemma}
\[A_{n}^{-1} =
 \prod_{i=0}^{n-1}a_{ii}^{-2} \] is given asymptotically by
 \[\exp\left\{d_{1}n^{2}\log n +d_{2}n^{2}+d_{3}n\log n
 +d_{4}n+d_{5}\log n +d_{6}+o(1)\right\}\]
 where
 $$d_{1}=1,\,\,\,d_{2}= -3/2,\,\,\, d_{3}=\nu,\,\,\,d_{4}=-\nu
 +\log 2\pi,\,\,\,d_{5}=\nu^{2}/2 -1/6$$
 $$ d_{6}= 4/3\log G(1/2) +(1/3 +\nu/2)\log \pi + (\nu/2 -1/18)\log2 -
\log G(1+\nu) $$
\end{lemma}
\proof
It is well known that the Barnes G function satisfies the property
\cite{WW}
$$G(1+z) = \Gamma(z)G(z).$$ From this it is quite easy to see
that $A_{n}^{-1}$ can be written as
\[\frac{G(1+n)G(1+n+\nu)}{G(1+\nu)}.\]
The asymptotics of the Barnes function are computed in \cite{WW}  and
since
$G(1+a+n)$ is asymptotic to
\[n^{(n+a)^{2}/2-1/12}e^{-3/4n^{2}-an}(2\pi)^{(n+a)/2}G^{2/3}(1/2)\pi^{1/6}
2^{-1/36}\]
we can directly apply this formula with $a=0$ and $a=\nu$ to obtain
the desired result.

This last result shows then that the asymptotics of the Hankel
matrices can be reduced to those of the matrices ${\cal H}_{n}(u) .$ We
compute these asymptotics by replacing the determinant with a
Fredholm determinant. To do this we must first consider some estimates
for the various kernels.

\section{Hilbert-Schmidt and trace norm estimates}

We proceed to write the determinant of ${\cal H}_{n}(u)$ as a Fredholm
determinant. It is clear that if $U(x)$ is a bounded function then the
matrix  ${\cal H}_{n}(u)$ can be realized as $P_{n} M_{U} P_{n}$ where
$M_{U}$ is multiplication by $U$ and $P_{n}$ is the projection onto
the space spanned by the first $n$ Laguerre functions. This is
clearly a bounded, finite rank operator defined on $L_{2}(0,\iy).$
Thus $$\det{\cal H}_{n}(u) =\det(I + K_{n})$$
is defined where $K_{n}$ is the integral operator whose kernel is given
by
\[ \sum_{i=0}^{n-1} L_{i}^{\nu}(x) L_{i}^{\nu}(y) (U(y) -1) .\]

Let us write a square root of the function $U -1$ as $V.$
Then
\[\det (I+K_{n}) =\det(I+P_{n}M_{V}M_VP_{n})=\det (I+ M_{V}P_{n}M_{V}).
\]
This last equality uses the general fact that $$\det(I+AB) =
\det(I+BA)$$ for general operators defined on a Hilbert space.
At this point we have replaced the kernel $K_{n}$ with
\[ \sum_{i=0}^{n-1} V(x) L_{i}^{\nu}(x) L_{i}^{\nu}(y) V(y)  .\]
Our primary goal in the paper is to replace this last kernel, which
we will commonly refer to as the Laguerre kernel, with a
more familiar one namely the compressed Bessel kernel, also defined
on $L_{2}(0,\iy),$
\bq
 V(x)
\frac{J_{\nu}(2(nx)^{1/2})\sqrt{ny}J_{\nu}^{'}(2(ny)^{1/2})
-J_{\nu}(2(ny)^{1/2})\sqrt{nx}J_{\nu}^{'}(2(nx)^{1/2})}{x-y}
V(y).\label{bes}
\eq
We use the term ``compressed'' since the above without the factors
$V(x)$ and $V(y)$ is generally called the
Bessel kernel.

We will use the general fact that if we have families of trace
class operators $A_{n}$ and
 $B_{n}$ (thinking of the Laguerre kernel as $A_{n}$
and the compressed Bessel kernel as $B_{n}$)  that depend on a parameter
$n$ such that
the Hilbert Schmidt
 norm $|| A_{n}-B_{n}||_{2} = o(1)$ and $\tr
 (A_{n}-B_{n})(I+B_{n})^{-1} =
 o(1)$ then
 \[\det(I+A_{n})/\det(I+B_{n})\,\,\ra 1\]
  as $n \ra \iy.$ The proof of this fact uses the idea of
 generalized determinants found in \cite{GK}, and also
 requires uniformity for the norms of the inverses of the
 operators $I+B_{n}.$ This will be made more precise later.

 Using the standard identity for the sum of the products of Laguerre
 functions allows us to write the Laguerre kernel
 \[ \sum_{i=0}^{n-1}V(x) L_{i}^{\nu}(x) L_{i}^{\nu}(y) V(y) \]
  as
\[
\sqrt{n(n+\nu)}V(x)\frac{L_{n-1}^{\nu}(x)L_{n}^{\nu}(y)-L_{n-1}^{\nu}(y)L_{n}^{\nu}(x)}
{x-y}V(y).\]
Notice that the above kernel has somewhat the same form as the
compressed Bessel kernel. It turns out that both the kernel above and
the Bessel kernel have an integral representation that make the
Hilbert-Schmidt computations simpler. The next lemma shows how this
is done for in the Laguerre case.

\begin{lemma}We have
\[{L_{n-1}^{\nu}(x)L_{n}^{\nu}(y)-L_{n-1}^{\nu}(y)L_{n}^{\nu}(x)\ov x-y}
={1\ov2}\int_0^1[L_{n-1}^{\nu}(tx)L_{n}^{\nu}(ty)+
L_{n-1}^{\nu}(ty)L_{n}^{\nu}(tx)]\,dt.\]
\end{lemma}
\proof Call the left-hand side $\Phi(x,y)$ and the right-hand
side $\Psi(x,y)$.
It follows from the differentiation formulas for Laguerre polynomials
that there is a
differentiation formula
\[x{d\ov
dx}\twoone{L_{n-1}^{\nu}(x)}{L_n^{\nu}(x)}=\twotwo{A}{B}{-C}{-A}\twoone
{L_{n-1}^{\nu}(x)}{L_n^{\nu}(x)},\]
where
\[A(x)={1\ov2} x -{\nu\ov2}-N, \ \ \
B(x)=C(x)=\sqrt{N(N+\nu)} \, . \]
An easy computation using this shows that $\Phi(x,y)$ satisfies
\[\left(x{\pl\ov\pl x}+y{\pl\ov\pl y}+1\right)\Phi(x,y)={1\ov2}
[L_{n-1}^{\nu}(x)L_{n}^{\nu}(y)+L_{n-1}^{\nu}(y)L_{n}^{\nu}(x)].\]
On the other hand, if we temporarily assume that $\nu>0$ and
differentiate under the
integral sign and then integrate by parts we find that
\[\left(x{\pl\ov\pl x}+y{\pl\ov\pl y}+1\right )\Psi(x,y)={1\ov2}
[L_{n-1}^{\nu}(x)L_{n}^{\nu}(y)+L_{n-1}^{\nu}(y)L_{n}^{\nu}(x)].\]
Hence
\[\left(x{\pl\ov\pl x}+y{\pl\ov\pl y}+1\right)(\Phi(x,y)-\Psi(x,y))=0.\]
{}From this it follows that $\Phi(x,y)-\Psi(x,y)$ is of the form
$\varphi(x/y)\,e^{-(x+y)/2}$
for some function $\varphi$ of one variable. Still assuming $\nu>0$,
both
$\Phi(x,y)$ and $\Psi(x,y)$ tend to 0 as $x$ and $y$ tend to 0
independently.
This shows that $\varphi=0$ and so the identity is established when
$\nu>0$. Since both
sides of the identity are analytic functions of $\nu$ for $\Re\,\nu>-1$
the identity holds
generally.

Thus our kernel may be rewritten
\be{1\ov2}\sqrt{n(n+\nu)}\,V(x)V(y)\,
\int_0^1[L_{n-1}^{\nu}(tx)L_{n}^{\nu}(ty)+L_{n-1}^{\nu}(ty)L_{n}^{\nu}(tx)]\,dt.
\label{kernel1}\ee
We assume from now on that $|\nu|\ge 1/2$.
The compressed Bessel kernel appearing in (\ref{bes}) also has the
following
well-known integral representation
\bq
nV(x)V(y)\int_{0}^{1}J_{\nu}(2\sqrt{nxt})J_{\nu}(2\sqrt{nyt})dt
\label{kernel2}
\eq
and the rest of the section is devoted to replacing the Laguerre
functions by the appropriate Bessel functions in (\ref{kernel1})
and showing that this leads to a small Hilbert-Schmidt error and then
showing that certain traces tend to zero. Notice that in the above
notation kernel (4) is the same as $A_{n}$ and kernel (5) is the same
as $B_{n}.$
We need only the following fact about the Laguerre functions.
\begin{lemma}\label{Lemma3.1}
Let $a$ be any real positive constant and suppose that $\nu \geq -1/2.$
Then as $n \ra \iy,$ the normalized Laguerre functions
satisfy
$$\max_{0<x \leq a}x^{1/4}L_{n}^{\nu}(x) =O(n^{-1/4}),$$
$$\max_{x \geq a}L_{n}^{\nu}(x) = O(n^{-1/4}).$$
We remark that the implied constants
in the estimates may depend on $a$ but not on $n$.
\end{lemma}
\proof These estimates follow easily from formula (7.6.9) and
Theorem 8.91.2 in \cite{Sz}.

\begin{lemma} We have
\be\l(x)-c_{n,\nu}\J(\st{x})={c_{n,\nu}\ov
6\sqrt\pi}N^{-3/4}x^{5/4}\sin(\st{x}-\al)+
O(N^{-5/4}(x^{-1/4}+x^3)),\label{estimate}\ee
where $N=n+(\nu+1)/2,
c_{n,\nu}=N^{-\nu/2}(\frac{\Gamma(n+\nu+1)}{\Gamma(n+1)})^{1/2}, \
\al=(\nu/2+1/4)\pi$ and the constant implied in
the O depends
only on $\nu$.
\end{lemma}
\proof Formula (8.64.3) of Szeg\"o reads, in current notation,
\[\l(x)-c_{n,\nu}\J(\st{x})\]
\[={\pi\ov 4 \sin\nu
\pi}\int_0^x[\J(\st{x})\Jm(\st{t})-\Jm(\st{x})\J(\st{t})]\,t\,\l(t)\,dt.\]
Using the Lemma \ref{Lemma3.1} estimate for $\l(t)$ we find that if
$Nx<1$ the
right side is at most
a constant times
\[(Nx)^{-|\nu|/2}\int_0^x(Nt)^{-|\nu|/2}\,t\,\l(t)\,dt
=O(N^{-|\nu|-{1\ov4}}x^{{7\ov4}-{|\nu|}})=O(N^{-2}).\]
(If $\nu$ is an integer this must be multiplied by $\log N$.)

So we assume $Nx>1$. The integral over $Nt<1$ is at most a constant
times
\[(Nx)^{-1/4}\int_0^{1/N}(Nt)^{-|\nu|/2}N^{-1/4}t^{{3\ov4}}\,dt
=O(N^{-{9\ov4}}x^{-1/4})\]
with, possibly, an extra $\log N$ factor. So we confine our integral now
to $Nt>1$.

If we replace the Bessel functions by their first-order asymptotics the
error is at
most a constant times
\[N^{-1}\int_0^x(x^{-1/2}+t^{-1/2})(xt)^{-1/4}N^{-1/4}(t+t^{3/4})\,dt
=O(N^{-5/4}(1+x)).\]
Therefore with this error we can replace the Bessel functions by their
first-order asymptotics,
obtaining, after using some some trigonometric identities,
\be{1\ov
4}N^{-1/2}x^{-1/4}\int_{1/N}^x\sin(\st{x}-\st{t})\,t^{3/4}\l(t)\,dt.
\label{lagint}\ee

Notice that this has the uniform estimate $O(N^{-3/4}(1+x^{3/2}))$, and
the
earlier errors
combined are $O(N^{-5/4}(x^{-1/4}+x)$. It follows that if in this
integral we replace $\l(t)$ by $c_{n,\nu}\J(\st{t})$ the error is
$O(N^{-5/4}(1+x^3))$.
Then if we replace $\J(\st{t})$ by its first-order asymptotics the error
is at most a
constant times
\[N^{-1/2}x^{-1/4}\int_{1/N}^xt^{3/4}\,(Nt)^{-3/4}\,dt=O(N^{-5/4}x^{3/4}).\]
Hence with the sum of the last-mentioned errors we may replace the
factor $\l(t)$
in (\ref{lagint}) by the
first-order asymptotics of $c_{n,\nu}\J(\st{t})$. Using a trigonometric
identity shows that this results in
\[{c_{n,\nu}\ov 8\sqrt\pi}N^{-3/4}x^{-1/4}\int_{1/N}^x[\sin(\st{x}-\al)-
\sin(\st{x}-4\sqrt{Nt}+\al)]\,t^{1/2}\,dt\]
\[={c_{n,\nu}\ov
6\sqrt\pi}N^{-3/4}x^{5/4}\sin(\st{x}-\al)+O(N^{-9/4}x^{-1/4})+
O(N^{-5/4}x^{3/4}).\]
Putting these things together gives the statement of the lemma.\sp

In what follows we use the notation $o_2(\,\cdot\,)$ or $O_2(\,\cdot\,)$
for a family
of operators whose Hilbert-Schmidt norm satisfies the corresponding $o$
or $O$
estimate, or for a kernel whose associated operator satisfies the
estimate.
Similarly, $o_1(\,\cdot\,)$ and $O_1(\,\cdot\,)$ refer to the trace
norm.
We also set $N'=n+(\nu-1)/2$.\sp

\begin{lemma}. Suppose $V$ is in $L^{\iy}$ and satisfies
$\int_{0}^{\iy}|V(x)|^{2}(x^{-1/2}+x^{6})dx < \iy$. Then the difference
between
the kernel (\ref{kernel1}) and
\[{1\ov2}\sqrt{n(n+\nu)}V(x)V(y)c_{n-1,\nu}
c_{n,\nu}\int_0^1[\J(\stp{tx})\J(\st{ty})+\J(\stp{ty})\J(\st{tx})]\,dt\]
is equal to

\noi (i) $O_1(N^{-1/2})$ plus a constant which is $O(1)$ times
\[N^{1/4}\int_0^1\left[(tx)^{5/4}\sin(\stp{tx}-\al)\J(\st{ty})
+(ty)^{5/4}\sin(\st{ty}-\al)\J(\stp{tx})\right]V(x)V(y)\,dt\]
plus a similar term with $x$ and $y$ interchanged;

\noi (ii) $O_2(N^{-1/4})$.
\end{lemma}
\proof The estimate of Lemma 3.3 holds with $n$ replaced by
$n-1$ if $N$ is
replaced by $N'$. Write
\[L_{n-1}^{\nu}(x)L_{n}^{\nu}(y)-c_{n-1,\nu}c_{n,\nu}\J(\stp{x})\J(\st{y})\]
\[=[L_{n-1}^{\nu}(x)-c_{n-1,\nu}\J(\stp{x})]\,\l(y)+L_{n-1}^{\nu}(x)[\l(y)-c_{n,\nu}\J(\st{x})]\]
\[+[c_{n-1,\nu}\J(\stp{x})-L_{n-1}^{\nu}(x)][\l(y)-c_{n,\nu}\J(\st{x})].\]
The contribution of the right side to the difference of the two kernels
is a constant
which is $O(1)$ times
\[N\int_0^1\left([L_{n-1}^{\nu}(tx)-c_{n-1,\nu}\J(\stp{tx})]\,\l(ty)\right.\]
\[+L_{n-1}^{\nu}(tx)[\l(ty)-c_{n,\nu}\J(\st{tx})]\]
\[+\left.[c_{n-1,\nu}\J(\stp{tx})-L_{n-1}^{\nu}(tx)][\l(ty)-c_{n,\nu}\J(\st{tx})]\right)
V(x)V(y)\,dt.\]
The integrand is a sum of three terms, each of which is the kernel (in
the $x,\,y$
variables) of a rank one operator. The trace norm of such an integral is
at most the
integral of the Hilbert-Schmidt norms.) Using this fact and Lemma 3.2
 we see that the contribution of the error term in (\ref{estimate})
to any of these terms is $O_1(N^{-1/2})$, as is the
contribution of the last term above. Then we see that replacing the two
Laguerre functions
by the corresponding Bessel functions leads to an even smaller error.
The error
term in Lemma 3.3 is seen also to contribute $O_1(N^{-1/2})$. Applying
this lemma,
and then doing everything with $x$ and $y$ interchanged we arrive at the
statement of
part (i).

For part (ii), we have to show that the integrals are $O_2(N^{-1/2})$.
It is easy to
see that with error $O_1(N^{-1/2})$ we may replace the Bessel functions
by their
first-order asymptotics, resulting in a constant which is $O(N^{-1/4})$
times
\[V(x)V(y)x^{5/4}y^{-1/4}
\int_0^1\sin(\stp{tx}-\al)\,\cos(\st{ty}-\al)\,t\,dt\]
plus similar expressions. A trigonometric identity and integration by
parts shows that
the integral is at most a constant times
\[(\max\,(1,\,|\stp{x}-\st{y}|))^{-1}+(\max\,(1,\,|\stp{x}+\st{y}|-2\al))^{-1},\]
and an easy exercise shows that, together with the outer factors, this
gives
$O_2(N^{-1/4})$.
Analogous argument applies to the other integrals and this
completes the proof of the lemma.\sp

\begin{lemma} The kernel (\ref{kernel1}) is equal to
\be
n\,V(x)V(y)\,\int_0^1\J(2\sqrt{ntx})\J(2\sqrt{nty})\,dt\label{bessel}\ee
plus an error $o_2(1)$.
\end{lemma}
\sp

\proof  Since the Laguerre kernel has Hilbert-Schmidt norm
$n^{1/2}$
(the operator is a rank $n$ projection) multiplying (\ref{kernel1}) by
constants which are
$1+O(n^{-1})$ produces an error $O_2(n^{-1/2})$. The constants we choose
are $n/(\sqrt{n(n+\nu)}c_{n-1,\nu} c_{n,\nu})$. It follows from this and
Lemma 3.4
that with error $o_2(1)$ we can replace (\ref{kernel1}) by
\be{n\ov2}V(x)V(y)\int_0^1[\J(\stp{tx})\J(\st{ty})+\J(\stp{ty})\J(\st{tx})]\,dt.
\label{bessel1}\ee

Now we show that if we replace $N$ and $N'$ by $n$ in this kernel the
error is $o_2(1)$.
Let's look  at the error incurred in the integral involving the
first summand when we replace $N'$ by $n$. It equals
\be\int_{N'}^n{dr\ov\sqrt
r}\int_0^1\sqrt{tx}J_{\nu}'(2\sqrt{rtx})\,\J(\st{ty})\,dt.
\label{rint}\ee

If $ntx<1$ and $nty>1$ the inner integral is at most a constant times
\[ n ^{-\mu/2-1/2}x^{-\mu/2}\int_0^{\min(1,\,1/nx)}t^{-\mu/2}dt,\]
where $\mu=\max(-\nu,\,0)$ as before. This is bounded by a constant
times
\[n^{-1/2}(nx)^{-\mu/2}\ {\rm if}\ nx<1,\ \ \ \ n^{-3/2}x^{-1}\ {\rm
if}\  nx>1.\]
This times $V(x)V(y)$ is $O_2(n^{-1})$ and so its eventual
contribution to the Hilbert-Schmidt norm (because of the external factor
$n$
and the fact that $r\sim n$) is $O(n^{-1/2})$.

If $ntx<1$ and $nty<1$ the inner integral is at most a constant times
\[n^{-\mu-1/2}x^{-\mu/2}y^{-\mu/2}\int_0^{\min(1,\,1/nx,\,1/ny)}t^{-\mu}dt.\]
By symmetry we may assume $y<x$. If $nx<1$ this is at most a constant
times the
outer factor which, when multiplied by $V(x)V(y)$, is $O_2(n^{-3/2})$.
If $nx>1$ the above is at most a constant times
\[n^{-3/2+\mu/2}x^{-1+\mu/2}y^{-\mu/2},\]
and this times $V(x)V(y)$ is $O_2(n^{-1})$. Thus the eventual
contribution of the portion of the integral where $ntx<1$ is
$O_2(n^{-1/2})$.

In the region $ntx>1$, if we replace the first Bessel function in
(\ref{rint}) by the
first term
of its asymptotic expansion it is easy to see we incur in the end an
error $O_2(n^{-1/4})$. After this replacement the inner integral becomes
a constant
times
\[n^{-1/4}x^{1/4}\int_{1/nx}^1t^{1/4}\cos(2\sqrt{rtx}-\al)\,\J(\st{ty})\,dt\]
\[=n^{-3/2}x^{1/4}\int_{1/x}^nt^{1/4}\cos(2\sqrt{rtx/n}-\al)\,\J(2\sqrt{ty})\,dt.\]
Taking account of the external factor of $n$ in our kernel, and the
$r$-integral, we
see that we want to show that the Hilbert-Schmidt norm of the kernel
\[V(x)V(y)\,n^{-1}\int_{1/x}^nt^{1/4}\cos(2\sqrt{rtx/n}-\al)\,\J(2\sqrt{ty})\,dt\]
tends to 0 as $n\ra\iy$. But from the asymptotics of the Bessel function
it is clear
that $n^{-1}$ times the integral is uniformly bounded by a constant
times $1+y^{-\mu/2}$ and tends to 0 whenever $x\ne y$. Hence the
dominated
convergence theorem tells us that the product is $o_2(1)$.

An analogous argument shows that replacing $N'$ by $n$ in the second
summand of
(\ref{bessel1}) and then replacing $N$ by $n$ in both summands leads to
an error $o_2(1)$.
This completes the proof.\sp

\noi{\bf Remark}. In the preceding lemmas our various kernels had the
factor $V(x)V(y)$.
It is easy to see from their proofs that the lemmas hold with this
factor replaced everywhere
by $V_1(x)V_2(y)$ as long as $V_1$ and $V_2$ are bounded and have
sufficiently rapid
decay at infinity. The next lemma is the first that will require some
smoothness.\sp

We shall denote by $J_n(x,y)$ the kernel $(\ref{kernel2})$ without the
external
factor $V(x)V(y)$, and by $J_n$ the corresponding operator. As before we
denote by $M_V$ multiplication by $V$ so that (\ref{kernel2}) is the
kernel
of the operator $M_V\,J_n\,M_V$.\sp

\begin{lemma} For any Schwartz function $W$ the commutator $[W,\,J_n]$
has Hilbert-Schmidt
norm which is bounded as $n\ra\iy$.
\end{lemma}
\proof Write the kernel of the commutator as in formula (3). It
has the form
\be
\frac{W(x) - W(y)}{x-y}(
J_{\nu}(2\sqrt{nx})\sqrt{ny}J_{\nu}^{'}(2\sqrt{ny}) )
\ee
plus a similar term with $x$ an $y$ interchanged. If $nx >1$ and $ny
>1$ then the product
$$ J_{\nu}(2\sqrt{nx})\sqrt{ny}J_{\nu}^{'}(2\sqrt{ny}) $$
is $O(x^{-1/4}y^{1/4}),$ and thus if $|x-y|$ is bounded away from
zero, or if we integrate over any bounded region
the Hilbert-Schmidt norm (the square root of the integral of the
square) is bounded.
If $|y-x| < 1,$ and assuming   $x, y > 1$ we see that $y/x$ is
bounded and thus the Hilbert-Schmidt norm can be estimated by the
square root of
\[ \int_{1}^{\iy}\int_{1}^{\iy}|\frac{W(x) - W(y)}{x-y}|^{2} \,dx\,dy.
\]
But this is known \cite{W2} to be bounded by
$\int_{-\iy}^{\iy}|x|\hat{W}(x)|^{2}dx.$
Now suppose that $nx <1$ and $ny >1.$ Then the product of Bessel
functions is $O(n^{-\mu+1/4}x^{-\mu/2}y^{1/4}).$ If $|x-y|$ is
bounded away from zero, then the resulting Hilbert-Schmidt norm is
$O(n^{-1/4}).$ If $|x-y| < 1,$ then it is also clear that the
Hilbert-Schmidt norm is $O(n^{-1/4})$. The other two cases
$xn >1, yn <1$  and  $xn <1, yn <1$ are handled
in the same fashion and are left to the reader.\sp

\nc{\sn}{\sqrt n} \nc{\Jt}{\tilde J_n}

\begin{lemma}For any bounded functions $W_1$ and $W_2$ we have
\[{\rm tr}\; (A_n-B_n)\,W_2\,J_n\,W_1\ra0\]
as $n\ra\iy$.
\end{lemma}
\proof For convenience all kernels $K(x,y)$ in this proof will be
replaced by
their unitary equivalents $2\sqrt{xy}K(x^2,\,y^2)$. We denote by
$L_n(x,y)$ this unitary
equivalent of the Laguerre kernel (\ref{kernel1}) without the external
$V$ factors
and by $J_n(x,y)$ here the unitary equivalent of the Bessel kernel
(\ref{kernel2})
without the $V$ factors. We also make the substitution $t\ra t^2$ in the
$t$ integrals.
Thus in our present notation
\be J_n(x,y)=4n\sqrt{xy}\int_0^1\J(2\sn xt)\,\J(2\sn
yt)\,tdt.\label{Jn}\ee
We also denote by $\Jt(x,y)$ the unitary equivalent of the first
displayed operator
of Lemma~3.4, without the $V$ factors. Thus
\[\Jt(x,y)=2\sqrt{n(n+\nu)}c_{n-1,\nu}c_{n,\nu}\]
\be\times\sqrt{xy}\int_0^1\J[(2\sqrt{N'}tx)\J(2\sqrt{N}ty)
+\J(2\sqrt{N'}ty)\J(2\sqrt{N}tx)]tdt.\label{Jt}\ee

If we set $V_i(x)=V(x^2)W_i(x^2)$ then
we see that our trace equals ${\rm tr}\,M_{V_1}(L_n-J_n)M_{V_2}J_n$. We
shall show
that this goes to 0
in two steps, showing first that ${\rm
tr}\,M_{V_1}(L_n-\Jt)M_{V_2}J_n\ra0$ and
then that ${\rm tr}\,M_{V_1}(\Jt-J_n)M_{V_2}\ra0$.

First, the asymptotics of the Bessel functions gives
\[\J(z)=\sqrt{2\ov\pi z}\cos(z-\al)+O\left({z^{-1/2}\ov<z>}\right),\]
where $<z>=(1+z^2)^{1/2}$. (This also uses $\nu\ge -1/2$.)
Hence
\[J_n(x,y)={8\ov\pi}n\sqrt{xy}\int_0^1\left[{\cos(2\sn tx-\al)\ov (2\sn
tx)^{1/2}}
+O\left({(\sn tx)^{-1/2}\ov<(\sn tx)>}\right)\right]\]
\[\times\left[{\cos(2\sn ty-\al)\ov (2\sn ty)^{1/2}}
+O\left({\sn ty)^{-1/2}\ov<(\sn ty)>}\right)\right]\,tdt\]
\[={8\ov\pi}n\int_0^1\cos(2\sn tx-\al)\,\cos(2\sn ty-\al)dt\]
\[+O\left({\sn\ov<\sn x>}+{\sn\ov<\sn x>}+{\sn\ov<\sn x>^{1/2}<\sn
y>^{1/2}}\right).\]
The last summand is at most a constant times the sum of the preceding
two. Using
this, a trigonometric identity and integrating we find that
\be J_n(x,y)=O\left({\sn\ov<\sn (x-y)>}+{\sn\ov<\sn (x+y)>}+{\sn\ov<\sn
x>}
+{\sn\ov<\sn y>}\right).\label{Jest}\ee

We consider first $M_{V_1}(L_n-\Jt)M_{V_2}J_n$. Lemma~3.4 tells us that
the kernel of $M_{V_1}(L_n-\Jt)M_{V_2}$ is $O_1(n^{-1/2})$ plus the
unitary
equivalent of the expression in part (i) with modified $V$ factors. (See
the remark following Lemma 3.5.) In the proof of part (ii) it
was stated that if we replace the Bessel functions in this expression
by their first order asymptotics the error is $O_1(n^{-1/2})$. (We shall
go through the
details for similar integrals below.) So we may replace $L_n-\Jt$ by a
constant which
is $O(1)$ times
\[\int_0^1[x^3\,\sin(2\sqrt{N'}tx-\al)\,\cos(2\sqrt{N}ty-\al)
+y^3\,\sin(2\sqrt{N'}ty\al)\,\cos(2\sqrt{N}tx-\al)]\,t^3dt.\]
(Recall the unitary equivalents we are using and the variable change
$t\ra t^2$.)
If in the integrals we made the replacements $N,\,N'\ra n$ we would
incur an error
$O((x^4+y^4)/\sn)$. Multiplying this by $V_1(x)V_2(y)$ times
(\ref{Jest}) and integrating
is easily seen to give $o(1)$. After these replacements the integral
becomes
what may be written
\[(x^3-y^3)\int_0^1\sin(2\sqrt{n}tx-\al)\,\cos(2\sqrt{n}ty-\al)\,t^3dt
+y^3\int_0^1\sin(2\sqrt{n}t(x+y)-2\al)\,t^3dt\]
\[=O\left({x^3-y^3\ov<\sn(x-y)>}+{x^3-y^3\ov<\sn(x+y)>}+{y^3\ov<\sn(x+y)>}\right)\]
\[=O\left({x^3-y^3\ov<\sn(x-y)>}+{x^3+y^3\ov<\sn(x+y)>}\right).\]
Let us see why if we multiply this by $V_1(x)V_2(y)$ times (\ref{Jest})
and
integrate we get $o(1)$.

First,
\[{x^3-y^3\ov<\sn(x-y>}{\sn\ov<\sn
(x-y)>}=O\left({x^2+y^2\ov<\sn(x-y)>}\right)\]
goes to zero pointwise and this times $V_1(x)V_2(y)$ is bounded by a
fixed $L^1$
function. Thus the integral of the product goes to zero. The term with
$<\sn(x+y)>$
instead of $<\sn(x-y)>$ is even smaller.

Next consider
\[{x^3-y^3\ov<\sn(x-y)>}{\sn\ov<\sn x>}.\]
We may ignore the factor $x^3-y^3$ since it may be incorporated into the
$V_i$.
After the substitutions $x\ra x/\sn,\ y\ra y/\sn$ the integral in
question becomes
\[{1\ov\sn}\int\int{|V_1(x/\sn)V_2(y/\sn)|\ov<x-y>\,<x>}dydx.\]
Schwarz's inequality shows that the $y$ integral is $O(n^{1/4})$ so our
double integral
is bounded by a constant times
\[n^{-1/4}\int_0^1dx+n^{-1/4}\int_1^{\iy}{|V_1(x/\sn)|\ov x}dx=n^{-1/4}+
n^{-1/4}\int_{1/\sn}^{\iy}{|V_1(x)|\ov x}dx=O(n^{-1/4}\log n).\]
Again the term with $<\sn(x+y)>$ instead of $<\sn(x-y)>$ is even
smaller.

Now we look at $M_{V_1}(\Jt-J_n)M_{V_2}J_n$. To find bounds for
$\Jt(x,y)-J_n(x,y)$
let us look first at the error incurred if in the first integral in
(\ref{Jt})
we replace $N'$ by $n$. The error in the integral together with the
external
factor $\sqrt{xy}$ equals
\[\int_{N'}^n{dr\ov\sqrt r}\int_0^1x^{3/2}y^{1/2}\J'(2\sqrt r tx)\J(2\sn
ty)\,t^2dt.\]
Using the asymptotics of $\J(z)$ and the fact that
\[\J'(z)=-\sqrt{2\ov\pi z}\sin(z-\al)+O(z^{-3/2}),\]
we can write the above as a constant times
\[\int_{N'}^n{dr\ov\sqrt r}\int_0^1x^{3/2}y^{1/2}\left[{\sin(2\sqrt r
tx-\al)
\ov (\sqrt r tx)^{1/2}}
+O((\sn tx)^{-3/2})\right]\]
\[\times\left[{\cos(2\sqrt N ty-\al)\ov (\sqrt N ty)^{1/2}}+
O\left({(\sn ty)^{-1/2}\ov<(\sn ty)>}\right)\right]t^2dt.\]
We estimate the trace norm of $V_1(x)V_2(y)$ times this by taking the
trace norm under
the integral signs. Since the integrand is, for fixed $r$ and $t$, a
function of $x$
times a function of $y$ its trace norm equals the product of the $L^2$
norms of its
factors. In multiplying out we will have main terms and error terms and
we must
estimate norms of all products. Thus we compute (in each line there will
be an integral
corresponding to a main term  and then and error term)
\[ \int{x^3|V_1(x)|^2\ov\sqrt r tx}dx=O(n^{-1/2}t^{-1}),\ \ \
\int{x^3|V_1(x)|^2\ov(\sn tx)^3}dx=O(n^{-3/2}t^{-3}),\]
\[ \int{y|V_2(y)|^2\ov\sqrt N ty}dy=O(n^{-1/2}t^{-1}),\ \ \
\int{y(\sqrt N ty)^{-1}|V_2(y)|^2\ov<\sqrt N ty>^2}dy=O(n^{-1}t^{-2}).\]
Combining $L^2$ norms we see that the trace norm of the contribution to
the integrand
of all but the product of the main terms is $O(n^{-3/4}t^{-2})$.
Integrating
over $t$ we are left with $O(n^{-3/4})$ and integrating over $r$ gives
$O(n^{-5/4})$.
If we combine this we the external factor in (\ref{Jt}) which is $O(n)$
we are left with
$O(n^{-1/4})$. Since the operator norms of the $J_n$ are bounded the
eventual contribution
to the trace of the product will be $O(n^{-1/4})$.

Thus we are left with the main term, which is
\[\int_{N'}^n(rN)^{-1/4}{dr\ov\sqrt r}\int_0^1x\sin(2\sqrt r tx-\al)
\,\cos(2\sqrt N ty-\al)\,tdt.\]
If in this we replaced $r$ and $N$ by $n$ everywhere in the integrand
the error would
be $O(n^{-3/2})$. If we multiply this by $V_1(x)V_2(y)$ and use the
estimate (\ref{Jest}) we
find by dominated convergence that the product has trace tending to
zero, even keeping in
mind the extra factor
$O(n)$ in (\ref{Jt}). So we may make these replacements, which results
in
\[(n-N')\,n\inv\int_0^1x\sin(2\sn tx-\al)\,\cos(2\sn ty-\al)\,tdt.\]

Now there is a second integral in $\Jt$, which is obtained from the
first by
interchanging $x$ and $y$. Interchanging and adding gives what can be
written
\[(n-N')\,n\inv\int_0^1[(x-y)\sin(2\sn tx-\al)
\,\cos(2\sn ty-\al)+y\,\sin(2\sn tx+2\sn ty-2\al)]\,tdt\]
\[=O\left({n^{-3/2}(x-y)\ov<\sn( x-y)>}+
{n^{-3/2}(x+y)\ov<\sn(x+y>}\right).\]
If we multiply by $V_1(x)V_2(y)$ and use the estimate (\ref{Jest}) we
find again that the
product has trace tending to zero, even keeping in mind the extra factor
$O(n)$.

Thus replacement of $N'$ by $n$ in (\ref{Jt}) leads to an eventual
error in the
trace of $o(1)$. Similarly so does then the replacement of $N$ by $n$.
Finally,
\[\sqrt{n(n+\nu)}c_{n-1,\nu}c_{n,\nu}=n\,(1+O(n\inv)),\]
so the eventual error in the trace upon replacing the constant by $n$ is
$O(n\inv)$
times what is obtained by multiplying $V_1(x)V_2(x)$ by the square of
(\ref{Jest})
and integrating. Dominated convergence shows this also to be $o(1)$.
This completes
the proof of the lemma.

\section{Completion of the proof}

Recall that $A_n$ is the Laguerre kernel (\ref{kernel1}) and $B_n$ is
the
compressed Bessel kernel (\ref{kernel2}), which we also denote in its
operator version as $M_V\,J_n\,M_V$. We shall show first that
\[\det\,(I+A_n)\sim\det\,(I+B_n)\]
as $n\ra\iy$. (It is clear that $A_{n}$ is a finite rank operator and so
it is trace class, and using the integral representation (\ref{bessel})
for the
compressed Bessel kernel and integrating over $t$ shows that $B_n$ is
also trace class.
Thus both determinants are defined.) This will follow from what we have
already done once we
know that the operators $I+B_n$ are uniformly invertible, which means
that
they are invertible for sufficiently large $n$ and the operator norms of
their
inverses are $O(1)$.\sp

\nc{\Vt}{\tilde V}

\begin{lemma} The operators $I+B_n$ are uniformly invertible and
\be(I+B_n)\inv=I-M_VJ_nM_{VU\inv}+O_2(1).\label{Binv}\ee
\end{lemma}
\proof We replace the operator by its unitary equivalent
$M_{\Vt}\,J_n\,M_{\Vt}$
where now $J_n$ is given by (\ref{Jn}), or equivalently
\[J_n(x,y)=\sqrt{xy}\int_0^{2\sn}\J(xt)\,\J(yt)\,tdt\]
and we set $\Vt(x)=V(x^2)$. If we set $H(x,y)=\sqrt{xy}\J(xy)$ with $H$
the corresponding
operator (the Hankel transform), and denote now by $P_n$ multiplication
by
the characteristic
function
of $(0,\,2\sn)$, then $J_n=H P_n H$ and so $I+B_n$ is unitarily
equivalent to
$I+M_{\Vt} H P_n H M_{\Vt}$. These operators will be uniformly
invertible if
$I+P_n H M_{\Vt^2} H P_n$ are.

Now it has recently been shown \cite{BE} that the operator
$H M_{\Vt^2} H$
is of the form $W(\Vt^2)+K$, where $W(\Vt^2)$ denotes the Wiener-Hopf
operator
with symbol $\Vt(x)^2=V(x^2)^2$ and $K$ is a compact operator on
$L^2(\bf R^+)$.
(Much less is needed for this than that $\Vt^2$ be a Schwartz function.)
Since
$1+V(x^2)^2=U(x^2)$ is nonzero and, being even, has zero winding
number it follows
from general facts about truncations of Wiener-Hopf operators that the
operators
$I+P_n W(\Vt^2) P_n$ are uniformly invertible. Then since $K$ is
compact
it follows that the $I+P_n (W(\Vt^2)+K) P_n$ will be uniformly
invertible
if the limiting operator $I+W(\Vt^2)+K=I+H M_{\Vt^2} H$ is invertible.
(For an exposition of the facts we used here see, for example, Chap.2 of
\cite{BS}.)
However, since
$H^2=I$ the inverse of $I+H\,M_{\Vt^2}\,H$ is easily seen to be
$I+H\,W\,H$ where
$W=-(1+\Vt^2)\inv\Vt^2$. This establishes the first statement of the
lemma.

For the second statement we apply Lemma 3.6 and use the facts
$J_n^2=J_n$ (which follows
from $H^2=I$) and that $V^2$ is a Schwartz function to see that for any
bounded
function $W$ we have
\[(I+B_n)(I+M_VJ_nM_W)=(I+M_VJ_nM_V)(I+M_VJ_nM_W)\]\[=I+M_VJ_nM_{V+V^2W+W}+O_2(1)
=I+M_VJ_nM_{V+UW}+O_2(1).\]
If we choose $W=-VU\inv$ and multiply both sides by $(I+B_n)\inv$ we
obtain the result.\sp

\begin{lemma} $\det\,(I+A_n)\sim\det\,(I+B_n)$ as $n\ra\iy$.
\end{lemma}
\proof If an operator $C$ is trace class then
\[\det(I+C) = \tilde{\det}(I+C)e^{-\tr C}\]
where $\tilde{\det}$ is the generalized determinant \cite{GK}.
(The generalized determinant is defined for any Hilbert-Schmidt
operator.)
Hence we can write
\[{\det\,(I+A_n)\ov\det\,(I+B_n)}=\det\,((I+A_n)(I+B_n)\inv)\]
\[=\det\,(I+(A_n-B_n)(I+B_n)\inv)=\tilde{\det}\,(I+(A_n-B_n)(I+B_n)\inv)\,
e^{-\tr(A_n-B_n)(I+B_n)\inv}.\]

It follows from Lemmas 3.4(ii) and 3.5 that $A_n-B_n\ra0$ in
Hilbert-Schmidt
norm, and therefore from the uniform invertibility of the $I+B_n$ that
the same
is true of $(A_n-B_n)(I+B_n)\inv$.
Therefore from the continuity of the generalized determinant in
Hilbert-Schmidt
norm that we conclude that the generalized determinant above has limit
1. Thus
it suffices to show that $\tr(A_n-B_n)(I+B_n)\inv\ra0$ as $n \ra \iy.$
Since $A_n-B_n\ra0$ in Hilbert-Schmidt the $O_2(1)$ term in (\ref{Binv})
 contributes  $o(1)$ to the trace of the product. By Lemma 3.7 the term
$M_VJ_nM_{VU\inv}$ in (\ref{Binv}) also contributes $o(1)$. That
$\tr(A_n-B_n)$ itself is $o(1)$ follows easily from arguments already
given---one
can check that at each stage the traces of the error operators tend to
zero.\sp

Finally, we can quote the main result of \cite{B} which gives the
asymptotics of $\det\,(I+B_n)$ or, more exactly the determinants of
their unitary equivalents. The formula
is
\[
 \det(I+B_{n}) \sim
 \exp  \{\frac{2n^{1/2}}{\pi}\int_{0}^{\infty}\log(U(x^{2})dx
 -\nu/2\log U(0) +\frac{1}{2\pi^{2}}\int_{0}^{\infty}xS(x)^{2}dx \}\]
where $S(x) = \int_{0}^{\iy}\cos(xy)\log(U(y^{2}) dy.$ This gives\sp

\begin{theorem} Suppose $U$ is nowhere zero and $U-1$ is a Schwartz
function.
 Then (\ref{asyformula}) holds.
 \end{theorem}

As mentioned in the introduction this result was computed
heuristically in \cite{BCW} using the Coulomb fluid approach. In the
same paper the analogous result was also obtained for weights
supported on the entire real line. These determinants involve Hermite
polynomials. It is highly likely that the results here (and
techniques) could also be
extended to that case.

\end{document}